%
%

\magnification=1200

\font\tpt=cmr10 at 12 pt
\font\fpt=cmr10 at 14 pt


\font\AAA=cmr14 at 12pt
\font\BBB=cmr14 at 8pt

\overfullrule=0in

\def\boxit#1{\hbox{\vrule
 \vtop{%
  \vbox{\hrule\kern 2pt %
     \hbox{\kern 2pt #1\kern 2pt}}%
   \kern 2pt \hrule }%
  \vrule}}

\def\Int{{\rm Int}}
\def\pn{\bbp^n}
\def\cn{\bbc^n}
\def\ss{\subset}

\def\dist{{\rm dist}}

\def\log{{\rm log}}

\def\arr{\longrightarrow}


\def\Theorem#1{\medskip\noindent {\AAA T\BBB HEOREM \rm #1.}}
\def\Prop#1{\medskip\noindent {\AAA P\BBB ROPOSITION \rm  #1.}}

\def\Note#1{\medskip\noindent {\AAA N\BBB OTE \rm  #1.}}
\def\Def#1{\medskip\noindent {\AAA D\BBB EFINITION \rm  #1.}}

\def\Ex#1{\medskip\noindent {\AAA E\BBB XAMPLE \rm    #1.}}

\def\pf{\medskip\noindent {\bf Proof.}\ }
\def\qed{\hfill  $\vrule width5pt height5pt depth0pt$}

   \def\cp{{\cal P}}   
   \def\co{{\cal O}}

\def\cl{{\cal L}}
\def\cp{{\cal P}}

\def\wt{\widetilde}
\def\wh{\widehat}

\def\and{\qquad {\rm and} \qquad}
\def\arr{\longrightarrow}

\def\bbc{{\bf C}}

\def\bbp{{\bf P}}

\def\a{\alpha}

\def\e{\epsilon}

\def\g{\gamma}

\def\z{\zeta}

\def\D{\Delta}

\centerline{\fpt   ANALYTIC DISKS AND THE PROJECTIVE HULL }
 
\vskip .2in
\centerline{\tpt Blaine Lawson$^*$ and John Wermer$^{**}$}
\vglue .9cm
\smallbreak\footnote{}{ $ {} \sp{ *}{\rm Partially}$  supported by
the N.S.F. } 
\smallbreak\footnote{}{ $ {} \sp{ **}{\rm Partially}$  supported by
the Institute Mittag-Leffler }

\centerline{\bf ABSTRACT} \medskip
  \font\abstractfont=cmr10 at 10 pt
Let $X$ be a complex manifold and  $\gamma$  a simple closed curve in                                   
$X$.  We address the question: What conditions on $\gamma$
insure the existence of a 1-dimensional complex subvariety $\Sigma$ with 
boundary $\gamma$ in $X$.                                                                                                                                                  
When  $X = \cn$, an answer to this question involves the                                           
polynomial hull of gamma.  When $X=\pn$, complex projective space,                                    
the projective hull $\wh\gamma$ of $\gamma$, a generalization of polynomial hull,     comes into play. One always has  $\Sigma\subseteq\wh\gamma$,  and for analytic $\gamma$
they conjecturally coincide.                            
 {{\parindent= .6in
\narrower\abstractfont \noindent

}}

  In  this paper we establish an approximate analogue of this  idea which holds
 without the analyticity of $\g$.
 We characterize points in $\wh\gamma$ as those  which
  lie on a sequence of analytic disks whose boundaries converge down
 to $\gamma$. This is in the spirit of work of Poletsky and  of   Larusson-Sigurdsson,
 whose results are essential here.
{{\parindent= .6in
\narrower\abstractfont \noindent

}}

 The results are applied to construct a remarkable example   
 of a closed curve $\gamma\ss\bbp^2$,   which is real analytic at all but one point, 
 and for which the closure of $\wh\gamma$     is $W\cup L$ where 
 $L$ is a projective line and $W$ is an analytic (non-algebraic) subvariety
 of $\bbp^2-L$. Furthermore, $\wh \g$ itself is   the union of $W$ with
 two points on $L$.                                                                                                      
{{\parindent= .6in
\narrower\abstractfont \noindent

}}

\vskip .3in

\centerline
{\bf   Introduction.}  
\medskip

Let $\gamma$ be a simple closed real curve in a complex manifold $X$.
Consider the problem of finding conditions which
guarantee  that  $\g$ forms the boundary of a complex analytic subvariety in $X$. 
When $X$ is $\cn$ (or, more generally, Stein),  there is 
a  solution  [W] which involves the polynomial hull of $\g$.
When $X$ is  $\pn$ (or, more generally, projective), there is a  notion of the projective
hull  of $\g$, denoted $\wh\g$, which is related to the polynomial hull and has 
the following property. If
$$
f:\Sigma \ \to\ X
$$
is a  map of a compact Riemann surface with boundary, which is holomorphic
on $\Int \Sigma$ and continuous up to the boundary with $$f(\partial \Sigma)=\g,$$
then 
$$
f(\Sigma)\ \subseteq \ \wh\g.
$$
(See [HL$_2$].)  
It is conjectured that if $X=\pn$ and 
$\g$ is real analytic, then either   $\wh\g=\g$
or $\wh \g$ is a 1-dimensional analytic subvariety   with boundary $\g$
or $\wh\g$ is an algebraic curve containing $\g$.
Real analyticity is fundamental for this conjecture.
It fails  for any $C^\infty$-curve $\g$ which is not pluripolar  [HL$_2$, Cor. 4.4].   
However, the conjecture has been established for real analytic curves
which are ``stable'' by using arguments of E. Bishop (see [HLW]).

In this paper we show that for any  closed curve $\g\ss\pn$,
all points $z_0$ in the projective hull $\wh\g$ are essentially characterized
by the fundamental property discussed above. We show that if $z_0\in\wh\g$, then
for any sequence $\e_j\searrow 0$, there exists a sequence of holomorphic maps of the unit disk
$$
f_j : \D\ \to\ \pn
$$
with 
$$
f_j(0)\ =\ z_0
$$
and
$$
\dist \left( f_j(\partial \D), \g\right)  \ \leq \ \e_j.
$$
In fact,  from any such sequence we extract limiting holomorphic maps
$f:\D\to \pn$, possibly constant $\equiv z_0$, and show    in Theorem 2 that $$f(\D)\ss\wh\g.$$

If we choose an affine coordinate chart $\cn\ss\pn$ containing $\g$, then the sequence
of maps  $f_j$ can be chosen so that the poles $\z^j_1,..., \z^j_{N_j}$ of $f_j$, viewed as
a $\cn$-valued function, are simple and satisfy
$
\sum_k \log|\z^j_k| \geq -M
$
for some constant $M$ independent of $j$.  In this form {\bf the existence of the sequence 
 is necessary and sufficient
for $z_0$ to lie in $\wh\g$.}  This is Theorem 1 below.
%
%

Theorem 2 is applied to produce interesting examples related to the conjecture
above. We construct a simple   closed curve $\g\ss   \bbp^2$, 
which is   real analytic at all but one point $p$, and has the following  properties.
There exists a projective line $L$ through $p$,   a point $q\in L$, and a proper  complex analytic subvariety $W\ss \bbp^2-L$   such that
$$
\wh\g \ =\ W \cup\{p,q\}  \and \overline{\left\{\wh\g\right\}}\ =\ W\cup L.
$$
Thus, if real analyticity  breaks down at a  single point, the conjecture fails.
Moreover,  if this happens, then $\wh\g$ may not be closed.

Related examples with other interesting properties are also given

One of us (John Wermer) wishes to acknowledge the contributions to                               
the genesis of this work he received from what he learned at                                       
Mittag-Leffler during the spring 2008 conference on Several Complex Variables.

\vfill\eject

\centerline
{\bf   The Main Theorems.}
\medskip

Consider a connected closed curve $\g$ lying in  complex projective $n$-space  $\pn$.
Write $\pn = \cn \cup H$, where $H$ denotes the hyperplane at infinity,
and assume that $\g\ss\cn$.  We  wish to characterize those points $z_0\in \cn$
which lie in  the projective hull  $\wh\g$ of $\g$  in terms of analytic disks.
We shall make use of the work of L\'arusson and Sigurdsson    [LS$_1$], who
  give a formula for the Siciak-Zahariuta extremal function $V_X$ of a connected open subset
$X$ of $\cn$.

For each $r>0$ let $K_r$ denote the open tube of radius $r$ around $\g$, i.e.,
$$
K_r\ =\ \{z \in\cn : {\rm dist}(z, \g)<r\}
$$
We fix a point $z_0\in \cn-\g$. Thus $z_0\notin K_r$ for $r$ sufficiently small.

Let $\{f_r\}_{r}$ be a family of analytic maps of the unit disk $\D$ into $\pn$,
indexed by numbers $r>0$ converging to zero. We consider
the following four conditions on the family  $\{f_r\}_{r}$:  

\medskip\noindent 
For all $r$
\medskip

(i)\ \ $f_r(\partial \D)\ \ss\ K_r$,

\medskip

(ii)\ \ $f_r(0)\ =\ z_0$,

\medskip

(iii)\ \ There exists a number $M>0$ such that if $\z^{(r)}_1,...,\z^{(r)}_{N_r}$ are the poles of 
$f_r$ in $\D$

\qquad (i.e., $f_r(\z^{(r)}_j) \in H$ for $j=1,..., N_r$), then we have
$$
\sum_{j=1}^{N_r} \log | \z^{(r)}_j |  \ \geq \ -M.
$$
\medskip

(iv)\ \ The poles are simple, that is, the  $\cn$-valued function $(\z-\z^{(r)}_j)f_r(\z)$ is holomorphic

\qquad\ in a neighborhood of $\z^{(r)}_j$ for each $j$.

\Theorem{1}  {\sl The point $z_0$ lies in $\wh\g$ if and only if there exists a family $\{f_r\}_r$  of analytic maps satisfying (i) -- (iv).}

\medskip\noindent
 {\bf Proof. (Sufficiency)}.
 Assume there exists a family  $\{f_r\}_r$ which satisfies (i) --  (iv).
 For given $r$, let $\z^{(r)}_1,...,\z^{(r)}_{N_r}$ be the poles of $f_r$.
 Put
 $$
 B_r(\z) \ \equiv \ \prod_{j=1}^{N_r}\left ({ \z-\z^{(r)}_j \over 1-\overline{\z}^{(r)}_j \z}    \right)
 \qquad{\rm for\ \ } \z\in\D.
 $$

\medskip\noindent
{\bf Claim:}  
$$
|B_r(0)|\ \geq \ e^{-M} \quad{\rm for\ all\ }r.
$$
\medskip\noindent
{\bf Proof of Claim.}  Since $B_r(0) = \prod_{j=1}^{N_r} (-\z^{(r)}_j )$ we have
$$
\log |B_r(0)| \ =\ \sum_{j=1}^{N_r} \log |\z^{(r)}_j |.
$$
By (iii) the right hand side is $\geq -M$.  Hence $|B_r(0)|\ \geq \ e^{-M} $ as claimed.\medskip

For each $r$ we now define a function $G_r$ by setting
$$
G_r\ \equiv \ f_r \cdot B_r.
\eqno{(1)}
$$
Note that $G_r$ has no poles on $\D$ and $|B_r|$ is of unit length on $\partial \D$.
  Also since $f_r(\partial \D) \ss K_r$, there exists a constant $c_1$ such that 
  $|G_r(\z)|\leq c_1$ on $\partial \D$ for all $r$. Hence, by the maximum principle
  $|G_r(\z)|\leq c_1$ on $\D$  for all $r$.  Thus $\{G_r\}_r$ is a normal family on the interior of
  $\D$, and so there exists a sequence $\{G_{r_j}\}_{r_j}$, $r_j\to 0$ as $j\to\infty$,
  converging point-wise on $\Int \D$ to a holomorphic  function $G$.
Furthermore, $\{B_r\}_r$ is a normal family on the interior of
  $\D$, and we may assume without loss of generality that 
$B_{r_j}\ \to B$, a holomorphic function on $\Int \D$.
The functions $B$ and $G$ lie in $H^\infty(\D)$.
Moreover, by Claim 1 we have $|B(0)|\geq e^{-M}$, so $B$ is not identically zero.

We put $Z$ = the zero set of $B$ on $\Int \D$.  Then $Z$ is a countable discrete subset
of $\Int\D$.  

Fix a point $a\in \Int \D \setminus Z$.  For each $r$, $f_r(a) = G_r(a)/B_r(a)$ and as 
$r_j\to0$, we have $B_{r_j}(a) \to B(a)$. By choice of $a$, $B(a)\neq 0$, so $\lim_{j\to\infty}f_{r_j}(a)$
exists and equals $G(a)/B(a)$.  We define
$$
f\ \equiv\ {G\over B}.
$$
Then $f$ is holomorphic on $\Int\D \setminus Z$ and has possible poles at the points of $Z$.
Since $f_r(0) = z_0\ \forall\ r$, we have
$$
f(0) = z_0.
$$
This brings us to the following.

\Theorem{2} {\sl For each $f$ constructed as above, we have
$$
f(\Int\D) \ \ss\ \wh\g.
$$
In particular, $f(0) = z_0\in \wh\g$.}

\medskip
\noindent
{\bf Proof of Theorem 2.}   We begin with the following.

\Prop{1} {\sl For each point  $a\in \Int \D \setminus Z$ we have}
$$
f(a) \in\wh\g.
$$

\pf
Fix a polynomial $P$ on $\cn$ with degree  $d$  Assume that $\|P\|_\g \leq 1$.
Then for all $r$ sufficiently small, $\|P\|_{K_r} \leq 2$

Fix such an $r$ and let $\z^{(r)}_1,...,\z^{(r)}_{N_r}$ be as above.  The map
$\z \mapsto P(f_r(\z))$ is holomorphic on $ \Int \D \setminus \bigcup_{j=1}^{N_r} \{ \z^{(r)}_j\}$.
Write $P$ as 
$$
P(w) \ =\ \sum_{|\a|\leq d} c_\a w_1^{\a_1} \cdots  w_n^{\a_n}.
$$
and $f_r(\z) = (w_1(\z),...,w_n(\z))$ so that 
$$
P(f_r(\z))\ =\ \sum_{|\a|\leq d} c_\a w_1^{\a_1}(\z) \cdots  w_n^{\a_n}(\z).
$$
Then by assumption (iv) for $\z$ near $\z^{(r)}_j$
$$
P(f_r(\z))\ =\ \left(  \z - \z^{(r)}_j \right)^{-d_j}\left ( a_j +b_j ( \z - \z^{(r)}_j) + c_j( \z - \z^{(r)}_j)^2+\cdots       \right)
$$
where $0\leq d_j\leq d$ and $a_j\neq 0$.  Hence,
$$
\log|P(f_r(\z))|\ =\ -d_j \cdot \log \left|  \z - \z^{(r)}_j \right| + h_j(\z)
\eqno{(2)}
$$
where $h_j$ is harmonic near $ \z^{(r)}_j$.   Also 
$$
\log \left| { \z - \z^{(r)}_j \over 1-   \overline{\z}^{(r)}_j \z}    \right|
\ =\ 
\log \left|  \z - \z^{(r)}_j \right| + k_j(\z)
$$
where $k_j$ is harmonic  on  $\D\setminus  \z^{(r)}_j$.  
We define
$$
\chi_r(\z) \ =\ \log |P(f_r(\z))| +d \cdot \log|B_r(\z)|\qquad {\rm for\ } \z\in\Int\D.
\eqno{(3)}
$$
On $\Int\D \setminus \bigcup_{j=1}^{N_r} \{  \z^{(r)}_j \}$ the function $\chi_r(\z)$ is
subharmonic.   

Fix $j=j_0$. For $\z$ near $\z^{(r)}_{j_0}$, (2) and (3) give
$$
\chi_r(\z) \ =\ -d_j \cdot \log \left|  \z - \z^{(r)}_{j_0} \right| + h_{j_0}(\z) + 
d \cdot \left \{\log \left| { \z - \z^{(r)}_{j_0} \over 1-   \overline{\z}^{(r)}_{j_0} \z}    \right|
+ \sum_{j\neq {j_0}}  \log \left| { \z - \z^{(r)}_j \over 1-   \overline{\z}^{(r)}_j \z}    \right |
   \right\}.
$$
Thus for $\z$ near $\z^{(r)}_{j_0}$,
$$
\chi_r(\z) \ =\  (d-d_{j_0})  \cdot \log \left|  \z - \z^{(r)}_{j_0} \right| + H_{j_0}(\z)
$$
where $H_{j_0}$ is subharmonic there.  Since $d \geq d_{j_0}$, the function
$\chi_r$ is subharmonic in a neighborhood of $\z^{(r)}_{j_0} $.

Since this holds for all $j_0$, $\chi_r$ is subharmonic on $\Int \D$. 
Also $\chi_r = \log |P(f_r)|$  on $\partial \D$ since $|B_r|=1$ there.
Therefore, for $\z\in\partial \D$ we have
$$
\chi_r(\z) \ =\ \log |P(f_r(\z))| \ \leq \log\, \| P \|_{K_r} \ \leq \ \log\, 2
$$
since $f_r(\partial \D) \ss K_r$ and $ \| P \|_{K_r} \leq 2$.  By the maximum principle 
for subharmonic functions on $\D$ we have that 
$$
\chi_r(\z) \ \leq \ \log\, 2\quad {\rm for\ } \z\in \D.
\eqno{(4)}
$$
Fix $a\in \Int \D \setminus Z$.  Then by (3), $\chi_r(a) \ =\ \log |P(f_r(a))| +d \cdot \log|B_r(a)|$, and so
$$
 \log |P(f_r(a))| +d \cdot \log|B_r(a)|   \ \leq\ \log\, 2.
$$
Letting $r\to0$, we get
$
 \log |P(f(a))| +d \cdot \log|B(a)|   \ \leq\ \log\, 2
$
and therefore
$$
 |P(f(a))| \ \leq \ 2\cdot \left|  {1\over B(a)} \right|^d.
\eqno{(5)}
$$
Now this holds for all polynomials $P$ with degree $\leq d$ and $\|P\|_\g \leq 1$.
Hence, $f(a) \in \wh \g$ and Proposition 1 is proved.\qed
\medskip

To complete the proof of Theorem 2 we must show that for each pole $\z_0$
of $f$ we have $f(\z_0)\in\wh\g$.  Fix a pole $\z_0$ and choose a small closed 
disk $D$ about $\z_0$ so that $f$ is regular on $D\setminus \{\z_0\}$. Let $\g_0 
= f(\partial D)$.  Then by [HL$_2$, Prop. 2.3] we have
$$
f(\z_0) \in \wh \g_0
$$
since $f(\z_0)$ lies on an analytic disk in $\pn$ with boundary $\g_0$.
This means that there is a constant $C_0>0$ such that for every section
$\cp \in H^0(\pn, \co(d))$ we have
$$
\|\cp(f(\z_0))\|\ \leq \ C_0^d  \|\cp\|_{\g_0}.
\eqno{(6)}
$$
Now in the affine chart $\cn\ss\pn$ each such $\cp$ corresponds to a 
polynomial $P$ of degree $\leq d$, and one has that for $z\in\cn$, 
$\|\cp\|_z = (1+\|z\|^2)^{-{d\over 2}}|P(z)|$.  It then follows from 
 (5) above that there is a constant $\kappa>0$
such that for all each $a\in \partial D$
$$
\|\cp(f(a))\|\ \leq \   \left|  {\kappa\over B(a)} \right|^d \|\cp\|_{\g}.
\eqno{(7)}
$$
Combining (6) and (7) gives a new constant $C>0$ such that
$$
\|\cp(f(\z_0))\|\ \leq \ C^d  \|\cp\|_{\g},
$$
which proves that $f(\z_0)\in \wh\g$ and establishes Theorem 2
and the sufficiency of conditions (i) -- (iv).\qed

\medskip\noindent
 {\bf Proof. (Necessity)}.  Assume $z_0\in \wh\g$.  
 We must provide a family
 $\{f_r\}_r$ of maps satisfying conditions (i) -- (iv).  The sets $K_r$ are defined as
 before.
%
%

\Def{1}  Let $E$ be a set in $\cn$ and denote by  $\cl_E$   the set  
functions $u$ in the Lelong class which satisfy   $u \leq 0$ on $E$.

\Def{2}  The Siciak-Zahariuta extremal function for $E$ is given by
$$
V_E(z) \ =\ \sup\{u(z) : u\in \cl_E\}.
$$

Since $z_0\in\wh\g$, we know from  [HL$_2$, p.6] that  there exists a constant $K$ such that 
$$
u(z_0)\ \leq\ K\qquad{\rm for\ all \ }\ u\in\cl_\g.
$$
 Fix $r$ and choose $u\in\cl_{K_r}$. Then $u$ is of Lelong class and $u\leq 0 $ on $K_r$.
 In particular $u\leq 0 $ on $\g$.   Hence $u\in \cl_\g$ and therefore $u(z_0)\leq K$.
From Definition 2 it follows that $V_{K_r}(z_0) \leq K$. Hence
$$
-V_{K_r}(z_0) \ \geq \  -K
\eqno{(8)}
$$

We now appeal to a result of L\'arusson and Sigurdsson on page 178 of [LS$_1$].
(Recall that $H=\pn\setminus\cn$ is the hyperplane at infinity.)

\Theorem {(L\'arusson - Sigurdsson)}
{\sl Let $X$ be a connected open subset of $\cn$.  Then for each $z\in\cn$,
$$
-V_X(z) \ =\ \sup_f \left( \sum_{f(\z)\in H} \log\,|\z|  \right)
$$
taken over all analytic maps $f:\D\to \pn$ with $f(\partial \D)\ss X$ and $f(0)= z_0$.}
\medskip

Because of (8) this theorem provides us with an analytic map $f_r:\D\to \pn$ with 
$f_r(\partial \D)\ss K_r$ and $f_r(0)=z_0$ such that 
$$
\sum_{f_r(\z)\in H} \log\,|\z| \ >\ -K-1.
\eqno{(9)}
$$
Putting $M=K+1$, we see that $f_r$ satisfies condition (iii).  (Note: $\{\z:f_r(\z)\in H\}
=\{\z^{(r)}_1,...,\z^{(r)}_{N_r}\}$.)
Standard transversality theory implies that for an open dense subset  of 
$G_{z_0}\equiv \{g\in {\rm PGL}(n+1,\bbc) :  g(z_0)=z_0\}$ the map $g\circ f_r$
is transversal to $H$, that is,  $g\circ f_r$ satisfies condition (iv).
(Since the group $G_{z_0}$ acts transitively on $\pn-\{z_0\}$, one can use, for example, 
 Sard's Theorem for families as  in [HL$_1$].)
Choosing $g$ sufficiently close to the identity we may assume that
$g\circ f_r(\partial \D)\ss K_r$ and that (9) holds with $f_r$ replaced by 
$g\circ f_r$. Choosing this approximation we see that $f_r'=g\circ f_r$ satisfies
all the conditions (i) -- (iv).

So we have constructed the desired family $\{f_r\}_r$, and necessity is established.
This completes the proof of Theorem 1.\qed

\Note{1}  The function $f$ appearing in Theorem 2 could be constant ($\equiv z_0$),
but if it is not, then we obtain a non-trivial analytic disk through $z_0$ which lies 
entirely in the projective hull $\wh\g$.

%
%

\medskip

\centerline{\bf The Examples.}
\medskip

\noindent
{\bf Example 1.}   We shall use Theorem 2 to construct  a closed curve $\g_0\ss\bbc^2$,
which is real analytic at all but one point, and whose projective hull contains a 
``large'' Riemann surface $\Sigma$.  In particular, $\Sigma$ will be the image
of a holomorphic map of the open disk $\Int \D$ which takes boundary values 
continuously (in fact, analytically) on $\g_0$ at all but one point. 
  We shall  then show
that the closure of  $\wh \g_0$ in $\bbp^2$
is exactly the union  of a projective line $L$ and a proper 
complex analytic (but not algebraic) subvariety 
 $W\ss \bbp^2\setminus L$ which extends $\Sigma$.

To produce $\Sigma$ we shall construct a family of holomorphic maps $f_n:\D\to \bbp^2$, $n=1,2,...$
satisfying conditions (i) -- (iv) for a sequence $r_n\searrow 0$ and with $z_0=(0,0)$.

To begin choose a sequence of  numbers $\{\e_j\}$  , $0<\e_j<1$, such that
$\sum_{j=1}^\infty \e_j < \infty$.  Put $a_j = 1-\e_j$, $j=1,2,...$.
Next choose a sequence of positive numbers $\{c_j\}$ such that
 $\sum_{j=1}^\infty  { c_j  \over \e_j} < \infty$.  For $n=1,2,... $ we put
 $$
 \omega_n (\z) \ =\    \sum_{j=1}^{n}{c_j\over \z-a_j}   + \sum_{j=1}^{n}{c_j\over a_j}.
 $$
Let  $f_n(\z) = (\z, \omega_n(\z))$ for $\z\in\D$.  Then $\{f_n\}$ is a sequence of holomorphic 
maps $\D\to \bbp^2$ such that $f_n$ has the poles $\z^{(n)}_j =a_j$, $j=1,..., n$. Put
$$
\omega(\z) \ =\ \sum_{j=1}^{\infty}{c_j\over \z-a_j}   + \sum_{j=1}^{\infty}{c_j\over a_j}
\qquad{\rm for\ } |\z|=1.
$$
For any $\z$ with $|\z|\geq 1$, we have $|\z-a_j|\geq 1-a_j$, so 
$$
\left|     {c_j\over \z-a_j}   \right|  \ \leq \ {c_j\over 1-a_j} \ =\ {c_j\over\e_j} 
$$
and so 
$$
\sum_{j=1}^{\infty} \left|     {c_j\over \z-a_j}   \right|  \ \leq \  \sum_{j=1}^{\infty}{c_j\over\e_j}, 
$$
and by our hypothesis the right hand side converges. Thus the series defining $\omega(\z)$ 
converges absolutely on $|\z|=1$. In fact it converges absolutely and uniformly for $|\z|\geq1$.

We define $\g_0$ to be the graph of the function $\omega$ over the curve $|\z|=1$ in
$\bbc^2$.  

Fix a point $\z$ with $|\z|=1$ and fix $n$.  The point $(\z,\omega(\z))$ lies on $\g_0$.
Hence,
$$
\eqalign
{
{\rm dist} (f_n(\z), \g_0)\ &\leq\ | \omega(\z) -   \omega_n(\z)   |    \cr
&=\    \left |   \left(   \sum_{j=1}^{\infty}{c_j\over \z-a_j}   + \sum_{j=1}^{\infty}{c_j\over a_j} \right)
-  \left(  \sum_{j=1}^{n}{c_j\over \z-a_j}   + \sum_{j=1}^{n}{c_j\over a_j}    \right)    \right |  \cr
&=\    \left |   \left(   \sum_{j=n+1}^{\infty}{c_j\over \z-a_j}   + \sum_{j=n+1}^{\infty}{c_j\over a_j} \right)\right|
\ \leq\  \sum_{j=n+1}^{\infty}{c_j\over 1-a_j}   + \sum_{j=n+1}^{\infty}{c_j\over a_j}.
}
$$

Fix $r$.  We recall that the set $K_r$ is the tube around $\g_0$ of radius $r$.  In view of the preceding,
$\dist(f_n(\partial \D), \g_0)$ becomes arbitrarily small for all $n$ large enough. So the 
sequence $\{f_n\}$ satisfies condition (i) for a suitable sequence of numbers $r_n\searrow 0$.

Next observe that $f_n(0) = (0, \omega_n(0))=(0,0)$ for all $n$. Hence, the sequence $\{f_n\}$
satisfies condition (ii) with $z_0=(0,0)$.

Finally, fix $n$ and note that  $\z^{(n)}_j = a_j$, $j=1,...,n$ are exactly the poles of the map
$f_n$. Now we have
$$
\log\, |\z^{(n)}_j |\ =\ \log\, a_j\ =\ \log(1-\e_j) \ \sim\ -\e_j,\ \ \ j=1,2,3,...
$$
so
$$
\sum_{j=1}^{n} \log\, |\z^{(n)}_j | \ \sim\  - \sum_{j=1}^{n} \e_j.
$$
Now $ \sum_{j=1}^{n} \e_j \leq  \sum_{j=1}^{\infty} \e_j \equiv M <\infty$ for all
$n$, and so
$
- \sum_{j=1}^{n} \e_j \geq -M
$
for all $n$.  Thus, for some $M'$ we have
$$
\sum_{j=1}^{n} \log\, |\z^{(n)}_j | \ \geq\ -M'\qquad{\rm for\ all\ \ } n,
$$
and the sequence satisfies condition (iii). Condition (iv) is straightforward to
verify, and we are done with the construction.

Fix a point $\z$ in $\D\setminus \bigcup_{j=1}^{\infty}a_j$.  Then $f_n(\z) = (\z,\omega_n(\z))$ and as 
$n\to\infty$,
$$
f_n(\z) \ \to\ (\z,\omega(\z)),
$$
where
$$
\omega(\z) \ 
=\ \sum_{j=1}^{\infty} { c_j \over \z-a_j  }   + \sum_{j=1}^{\infty}{c_j\over a_j}.
\eqno{(10)}
$$
It is easily verified that this series converges uniformly  in $\z$ on compact subsets of $\Int \D
\setminus \bigcup_{j=1}^{\infty} a_j$. In fact it converges uniformly on compact subsets
of the domain $\bbc \setminus \{1\}\cup \bigcup_{j=1}^{\infty} a_j$.

Consider the meromorphic map $f(\z)= (\z,\omega(\z))$ on $\Int \D$. It follows from Theorem 2 that
$$
\Sigma \ \equiv\ f(\Int\D) \ \ss\ \wh \g_0.
$$
This includes all points on the graph of $\omega$ over $\Int \D\setminus \bigcup_{j=1}^{\infty}a_j$.

\medskip\noindent
{\bf NB.  }  The meromorphic map $f$ is in fact a holomorphic map $f: \Int \D\to \bbp^2$.
However, its image passes {\bf infinitely often} through the point $\ell \in H\cong \bbp^{1}$ 
corresponding to the ``vertical'' line in $\bbc^2$.  In particular, the image of $f$ is not an
analytic subvariety at that point.
\medskip

We   observed above  that the series for  
$
\omega(\z) 
$
given in (10) converges uniformly on the set $\bbc \setminus \Int \D$.  Moreover,  its graph extends across infinity to give a regularly embedded disk $\Sigma^-$ in $\bbp^2$ with boundary $\g_0$, taken
from the ``outside''. Thus by [HL$_2$, Prop. 2.3] we have $\Sigma^-\ss\wh\g_0$.

We now denote by $L\ss\bbp^2$ the projective line determined by $\z=1$.
Let $W$ be the closure in $\bbp^2-L$ of the graph of $\omega$. 
Note that 
$$W=\Sigma\cup \Sigma^- \setminus\{\ell, p\}$$
where $\ell$ is the common polar point referred to above and $p=(1,\omega(1))$.
Note that  $W$ is
a complex analytic subvariety of dimension 1 in $\bbp^2-L$.  We have proved that
$W\subseteq\wh \g_0$.

\Theorem {3}  {\sl For appropriate choices of the sequences $\{\e_j\}$ and $\{c_j\}$ one has that}
$$
\wh \g_0   \ =\ W \cup\{\ell,p\} \and  \overline{(\wh{\g_0})} \ =\ W\cup L
$$

\pf
We must show that points of $\bbp^2\setminus (W \cup\{\ell,p\})$ do not lie on $\wh\g_0$. 
The second assertion then follows from the Picard Theorem.

Consider the polynomial of degree $N+1$:
$$
P_N(z,w) \  \equiv \left(w-\kappa - \sum_{n=1}^{N} {c_n\over z-a_n} \right )\prod_{n=1}^{N}(z-a_n)
\eqno{(11)}
$$
where $\kappa = \sum_{n=1}^\infty c_n/a_n$.
Note that 
$$
P_N(z,\omega(z))\ =\ \left (\sum_{n=N+1}^{\infty} {c_n\over z-a_n} \right) \prod_{n=1}^{N}(z-a_n)
$$
For $|z|=1$ we have the estimate that $|z-a_n|\geq \e_n$. Hence we have
$$
\|P_N\|_{\g_0} \ =\ \sup_{|z|=1}\left|    \left (\sum_{n=N+1}^{\infty} {c_n\over z-a_n} \right) \prod_{n=1}^{N}(z-a_n)   \right| \ \leq\      \left (\sum_{n=N+1}^{\infty} {c_n\over  \e_n} \right) 2^N
\eqno{(12)}$$

Now choose $\{c_n\}, \{\e_n\}$ so that 
$$
\sum_{n=N+1}^{\infty} {c_n\over  \e_n} \ <\ \left({1\over N+1}\right)^{N+1}
\eqno{(13)}$$
For example set $\e_n = {1\over 2^n}$ and $c_n =  {1\over 2^n} {1\over 2^n} {1\over n^n}$

Now choose $z$ with $z\neq1$ and $z\neq a_n$ for any $n$. Pick any $w\neq \omega(z)$.
Consider equation (11).
The first factor on the RHS converges to $w-\omega(z)\neq 0$, and so its
$(N+1)^{\rm st}$ root converges to 1.  The second factor satisfies
$$
\left | \prod_{n=1}^{N}(z-a_n)  \right|^{1\over N +1} \ \arr \ |z-1|.
$$
(These assertions follow from the estimate: $|(a + \e)^{1\over n} -a^{1\over n}| \leq {c\e\over n}$
for $0<|\e| < a/2$ where   $c$ depends only on $a$  , and on our assumption that $\sum \e_k<\infty$.)
Thus we have that 
$$
\lim_{N\to\infty} |P_N(z,w)|^{1\over N+1}\ =\  |z-1|\ \neq \ 0.
$$
On the other hand, we have by (12) and (13)  that 
$$
\left\{   \|P_N\|_{\g_0} \right\}^{1\over N+1}\  \leq \ {2\over N+1}  \ \to\ 0\quad {\rm as\ \ } N\to \infty.
$$
It follows that there cannot exist a constant $C>0$ so that 
$$
|P_N(z,w)|^{1\over N+1} \ \leq \ C    \left\{   \|P_N\|_{\g_0} \right\}^{1\over N+1} \qquad {\rm for\ all\ \ } N,
$$
which means that  $(z,w)\notin \wh{\g_0}$ as claimed.

Suppose now that $z=a_n$ for some $n$.   In equation (11) move   $(z-a_n)$  over to the left 
factor  so that  factor becomes regular at $a_n$. Then we have
$$
P_N(a_n,w) \ =   \ c_n  \prod_{j\neq n}^{N}(a_n-a_j).
$$
Again we see that  
$$
\left | \prod_{j\neq n}^{N}(a_n-a_j)  \right|^{1\over N+1} \ \arr \ |a_n-1| \ \neq\ 0.
$$
and the same contradiction results.

Suppose now that $\z=1$ and $w\neq \omega(1)$. We now choose our sequence 
$\{c_n\}$ to converge even more rapidly  so that 
$$
\sum_{n=N+1}^{\infty} {c_n\over  \e_n} \ <\ \left({1\over N+1}\right)^{(N+1)(N+1)}
\eqno{(13)'}
$$
For example set $\e_n = {1\over 2^n}$ and $c_n =  {1\over 2^n} {1\over 2^n} {1\over n^{n^2}}$.
Then for large $N$ one has 
$$
P_N(1,w) \ \sim\ (w-\omega(1))\prod_{n=1}^{N} \e_n \ =\ (w-\omega(1))\prod_{n=1}^{N}  {1\over 2^n} 
\ =\   (w-\omega(1)) \left({1\over 2}\right)^{N(N+1)\over 2}
$$
Comparing with (13)$'$ as above shows that $(1,w) \notin \wh\g_0$ when $w-\omega(1)\neq 0$.

It now remains only to eliminate all  points on the line $H$ at infinity except for $\ell$ and $p$.
We begin with the following observation.  Let $\cp_N\in H^0(\bbp^2,\co{}(N+1))$
denote the holomorphic section corresponding to the polynomial $P_N$.
Then equations (12) and (13) imply that for some constant $K$ and all $N$
$$
\sup_{\g_0} \| \cp_N\| \ \leq \  \left({2K\over N+1}\right)^{N+1}
\eqno{(14)}
$$
where $\|\bullet \|$ denotes the standard metric in the line bundle $\co{}(N+1)$.
This equation (14) can be interpreted in any coordinate chart.

We make a change of coordinates as follows. First let  $s=z-1$ and set $  P_N'(s,w) = P_N(s+1,w)$.
We now pass to homogeneous coordinates $(t_0,s_0,w_0)$ where the corresponding
homogeneous polynomial  is 
$$
Q_N(t_0,s_0,w_0) \ \equiv \  t_0^{N+1}  P_N'\left({s_0\over t_0}, {w_0\over t_0}\right).
$$
Next we pass to the affine coordinate chart defined by setting $s_0=1$, or equivalently
by dividing by $s_0$.  This gives new coordinates $(t_1,w_1)$ where $t_1=t_0/s_0$ 
and $w_1=w_0/s_0$. Thus, the change of coordinates from the old chart  (where
$t_0=1$) is: $t_1=1/s$, $w_1=w/s$.

For simplicity of notation we relabel these new affine  coordinates as $(t,w)$.
In this  affine chart our polynomial is expressed in terms of $Q_N$ be setting
$s_0=1$, that is, the polynomial is now
$P_N''(t,w) = Q_N(t,1,w)$.   Calculation shows that
$$
P_N''(t,w)\ =\ \left( w-\kappa t - \sum_{n=1}^{N} {c_n t^2\over 1+\e_n t}\right) \prod_{n=1}^{N}
(1+\e_nt).
$$
Now in the affine $(t,w)$ coordinates the line $L$ has become the line 
at infinity, and the old line at infinity  $H$ corresponds to $\{t=0\}$.
The point $\ell$ lies at infinity on $H$ and the point $p$ corresponds to $(0,0)$.
Note that
$$
P_N''(0,w)\ =\ w \and \|\cp_N''(0,w)\| \ =\  \left( {1\over 1+|w|^2}\right)^{{N+1\over 2}} |w|
\eqno{(15)}
$$
Now if $(0,w) \in \wh\g_0$ for $w\neq0$, then there would be a constant $C>0$ such that
$$ 
\|\cp_N''(0,w)\|^{1\over N+1} \ \leq\  C \left(  \sup_{\g_0} \| \cp_N''\|  \right)^{1\over N+1}
$$
contradicting (14) and (15).\qed

 \medskip

\Ex{2}   We repeat the construction above with poles clustering at all points of $\partial \D$.
Put
$$
\wt\omega_n(\z)   \ \equiv\   \sum_{k=1}^{n}  \sum_{\ell=1}^{k} 
\left( {c_k  \over  \z-e^{2\pi i \ell\over k}a_k} \right)- \kappa_n
$$
where $\kappa_n$ is chosen so that $\wt\omega_n(0)=0$.
Let $a_k = 1-\e_k$ and choose $\e_k>0$ and $c_k>0$ so that
$\sum_k\e_k<\infty$ and $\sum_k{kc_k\over \e_k}<\infty$. We   now
proceed in exact analogy with Example 1. 
The limit  $\omega = \lim_n \omega_n$ converges absolutely on $\partial \D$
and its graph defines a curve $\g_\infty$ in $\bbc^2$. The same limit over $\Int \D$ defines
a meromorphic function whose graph lies in the projective hull   $\wh\g_\infty$ by 
Theorem 2.  This limit also exists at all points of $\bbc\setminus \D$ and gives
an exterior analytic disk contained in $\wh\g_\infty$. 

In this example the closure of $\wh\g_\infty$ contains $\partial\D\times \bbc$, a subset
if dimension 3.

Set 
$$
\omega(\z) \ =\ \sum_{n=1}^{\infty} {c_n\over \z-a_n} + \kappa \qquad{\rm where\ }
\kappa\ =\  \sum_{n=1}^{\infty} {c_n\over  a_n}
$$
We are considering the graph $W$ of $f(\z)=(\z,\omega(\z))$ for $\z\neq a_n$ any $n$.
Our curve $\g_0$ is just the graph of $\omega$ above $\partial \D$.
For rapidly converging $\{c_n\}$ the analogue of Theorem 3 will hold.

\vfill\eject

\centerline{\bf 
References.}\bigskip

\item{[HL$_1$]} F. R. Harvey and H. B. Lawson, Jr. {\sl On boundaries of complex analytic varieties, I},
Annals of  Math. {\bf 102} (1975),  223-290.

\smallskip

\item{[HL$_2$]} F. R. Harvey and H. B. Lawson, Jr. {\sl Projective hulls and the projective Gelfand transform},
Asian J. Math. {\bf 10} (2006), 607-646.
\smallskip

\item {[HLW]}  F. R. Harvey, H. B. Lawson, Jr. and J. Wermer, 
{\sl  The projective hull of certain curves in $\bbc^n$},\ 
(with  ), Stony Brook Preprint (2006). ArXiv:math.CV/0611482.


\smallskip
\item{[LS$_1$]}  F. Larusson and R. Sigurdsson,  {\sl The Siciak-Zahariuta extremal                       
function as the envelope of disc functionals}, Ann. Polon. Math. {\bf 86}                         
(2005), 177-192.

\smallskip
\item{[LS$_2$]}   F. L\'arusson and R. Sigurdsson, {\sl Siciak-Zahariuta extremal functions and polynomial
hulls}, Annales Polonici Math. (2007), 235-239.

\smallskip

\item{[P]}   E. A. Poletsky,  {\sl Holomorphic currents}, Indiana Univ.  Math. J. (1993),  85-144.

\smallskip

 \noindent
[W]   J. Wermer  {\sl    The hull of a curve in $\bbc^n$},    
Ann. of Math., {\bf  68}  (1958), 550-561.

 \end

Dear Blaine,                                                                                 
  here are my comments on the latest version.                                                
                                                                                             
  where we refer to [LS] on our page 5 (just below (8)), the reference                       
is incorrect. We have to add, to our list of references,:                                    
                                                                                             
[LS_{1}]  F. Larusson and R. Sigurdsson, The Siciak-Zahariuta extremal                       
function as the envelope of disc functionals" , Ann.Pol.Math. 86.2                           
(2005)                                                                                       
                                                                                             
The result we quote is on p.178 of that paper. The other one, I think,                       
should be left in our references, as [LS_{2}].                                               
                                                                                             
I have been reading the end  of our manuscript, your new results                             
regarding points on  H.  First, congratulations! Still , I have                              
difficulties with it.                                                                        
                                                                                             
The change of coordinates from (z,w) to (t,w) should be made explicit.                       
                                                                                             
On p. 10, line 6 from top,, it should be , on the right side of the equation,                
                                                                                             
||P_{N}"|| rather than ||P_{N}||.  If not , please explain.                                  
                                                                                             
       Regards,  John

Dear Blaine,

Regarding the Introduction, I propose:

To leave out the part starting with  "In this paper..., and ending with
"whose work [LS] is
essential here ", on p.2, since most of that material is covered in the
Abstract and on p.3.

Also, please leave out the words  "an important solution, which goes
back to one to one of the authors [W]."

Let's talk soon on the phone !

                Best regards,  John

Dear Blaine:

I like the latest version of the change of coordinates

I suggest the folowing tiny addition, to the Abstract, 2nd paragraph:

   "...holds without analyticity"  should be replaced by  "holds  
without analyticity of gamma".

   I suggest we now submit the manuscript to  Mittag-Leffler 2008  
Spring Preprint Series , and send it to Bo Berndtsson at
                 berndtsson@mittag-leffler.se,

together with the following sentence:

    Ths enclosed is joint work by the two of us, submitted to the  
Mittag-Leffler
Spring 2008 preprint series.
     We also want to put an Arxiv version to the ARxiv.  We hope that  
is OK with you.

                                 Regards, and then both our names.

        hoping to speak to you soon on the phone.  Greetings!  John

\end